\documentclass[11pt]{article}
\usepackage{graphicx}
\usepackage{amsmath,amsfonts,amssymb,amsthm}
\usepackage{bbm,mathrsfs, bm, mathtools}

\title{The polytopes in a Poisson hyperplane tessellation}
\author{Rolf Schneider}
\date{}
\newcommand{\Sd}{{\mathbb S}^{d-1}}
\newcommand{\R}{{\mathbb R}}

\newcommand{\K}{{\mathcal K}}

\newcommand{\bP}{{\mathbb P}}

\newcommand{\N}{{\mathbb N}}

\newcommand{\Ha}{\mathcal{H}}

\newcommand{\B}{\mathcal{B}}
\newcommand{\D}{{\rm d}}

\newcommand{\bE}{{\mathbb E}\,}

\newcommand{\bQ}{{\mathbb Q}}

\newtheorem{theorem}{Theorem}
\newtheorem{proposition}{Proposition}
\newtheorem{definition}{Definition}

\begin{document}
\maketitle

\begin{abstract}
For a stationary Poisson hyperplane tessellation $X$ in $\R^d$, whose directional distribution satisfies some mild conditions (which hold in the isotropic case, for example), it was recently shown that with probability one every combinatorial type of a simple $d$-polytope is realized infinitely often by the polytopes of $X$. This result is strengthened here: with probability one, every such combinatorial type appears among the polytopes of $X$ not only infinitely often, but with positive density. 
\\[2mm]

\end{abstract}

\section{Introduction}\label{sec1}

Imagine a system ${\sf H}$ of hyperplanes in Euclidean space $\R^d$ ($d\ge 2$) that induces a tessellation $T_{\sf H}$ of $\R^d$. This means that any bounded subset of $\R^d$ meets only finitely many hyperplanes of ${\sf H}$ and that the components of $\R^d\setminus \bigcup_{H\in {\sf H}} H$ are bounded. The closures of these components are then convex polytopes which cover $\R^d$ and have pairwise no common interior points. The set of these polytopes is denoted by $T_{\sf H}$. We impose the additional assumption that the hyperplanes of ${\sf H}$ are in general position; then each polytope of $T_{\sf H}$ is simple, that is, each of its vertices is contained in precisely $d$ facets. The polytopes appearing in $T_{\sf H}$ may be rather boring; they could, for example, all be parallelepipeds. However, if the hyperplanes of ${\sf H}$ have sufficiently many different directions, one can imagine that quite different shapes of polytopes appear in $T_{\sf H}$. Is it possible that every combinatorial type of a simple $d$-polytope is realized in $T_{\sf H}$? This can be achieved in a much stronger sense. 

In fact, suppose that $\widehat X$ is a stationary and isotropic Poisson hyperplane process in $\R^d$ (explanations are found in \cite{SW08}, for example). Its hyperplanes are almost surely in general position and induce a random tessellation of $\R^d$, denoted by $X$. The general character of the polytopes in $X$ was recently investigated in \cite{RS16}. For example, it was shown there that almost surely (a.s.) the translates of the polytopes in $X$ are dense in the space of convex bodies in $\R^d$ (with the Hausdorff metric). Another result was that a.s. the polytopes of $X$ realize every combinatorial type of a simple $d$-polytope infinitely often. In the following, we improve the latter result considerably, replacing `infinitely often' by `with positive density'. In the subsequent definition, $B_n$ is the ball in $\R^d$  with center at the origin and radius $n\in\N$, and $\lambda_d$ denotes Lebesgue measure in $\R^d$. Further, ${\mathbbm 1}_A$ is the indicator function of $A$.

\begin{definition}
Let $T$ be a tessellation of $\R^d$, and let $A$ be a translation invariant set of polytopes in $\R^d$. We say that $A$ {\bf appears in $T$ with density $\delta$} if
$$ \liminf_{n\to \infty} \frac{1}{\lambda_d(B_n)} \sum_{P\in T,\,P\subset B_n} {\mathbbm 1}_A(P)=\delta.$$
\end{definition}

With this definition, we prove below that in a Poisson hyperplane tessellation in $\R^d$ which is stationary and isotropic (that is, has a motion invariant distribution), almost surely every combinatorial type of a simple $d$-polytope appears with positive density. The actual result will, in fact, be more general: it is sufficient that the Poisson hyperplane tessellation is stationary and that its directional distribution, a measure on the unit sphere, is not zero on any nonempty open set and is zero on any great subsphere. The precise theorem is formulated in the next section.

\section{Explanations}\label{sec2}

We work in the $d$-dimensional Euclidean space $\R^d$ ($d\ge 2$) with its usual scalar product $\langle\cdot\,,\cdot\rangle$. By $\lambda_d$ we denote its Lebesgue measure, by $o$ its origin, by $B^d$ its unit ball (with $nB^d=:B_n$), and by $\Sd$ its unit sphere. The space of hyperplanes in $\R^d$, with its usual topology, is denoted by $\Ha$, and $\B(\Ha)$ is the $\sigma$-algebra of Borel sets in $\Ha$. Hyperplanes in $\R^d$ are often written in the form
$$ H(u,\tau)= \{x\in \R^d: \langle x,u\rangle \le\tau\}$$
with $u\in\Sd$ and $\tau\in \R$.

We assume that $\widehat X$ is a stationary Poisson hyperplane process in $\R^d$, thus, a Poisson point process in the space $\Ha$ of hyperplanes, with the property that its distribution is invariant under translations (we refer, e.g., to \cite{SW08} for more details). The {\em intensity measure} $\widehat\Theta$ of $\widehat X$ is defined by
$$ \widehat\Theta(A)=\bE \widehat X(A)\quad\mbox{for }A\in\B(\Ha).$$
Here $\bE$ denotes expectation, and we write $(\Omega,{\mathcal A},\bP)$ for the underlying probability space. It is assumed that $\widehat\Theta$ is locally finite and not identically zero. That $\widehat X$ is a Poisson process includes that
$$ \bP(\widehat X(A)=k)= e^{-\widehat \Theta(A)}\frac{\widehat \Theta(A)^k}{k!} \quad \mbox{for }k\in \N_0,$$
for any $A\in \B(\Ha)$ with $\widehat\Theta(A)<\infty$.

Since $\widehat X$ is stationary, the measure $\widehat\Theta$ has a decomposition
$$ \widehat\Theta(A) =\widehat \gamma \int_{\Sd} \int_{-\infty}^\infty {\mathbbm 1}_A(H(u,\tau))\,\D\tau\,\varphi(\D u)$$
for $A\in\B(\Ha)$ (see \cite{SW08}, Theorem 4.4.2 and (4.33)). The number $\widehat\gamma>0$ is the {\em intensity} of $\widehat X$, and $\varphi$ is a finite, even Borel measure on the unit sphere. It is called the {\em spherical directional distribution} of $\widehat X$. For any such measure $\varphi$ and any number $\widehat\gamma>0$, there exists a stationary Poisson hyperplane process in $\R^d$ with these data, and it is unique up to stochastic equivalence. 

The hyperplane process $\widehat X$ induces a random tessellation of $\R^d$, which we denote by $X$. As usual, a random tessellation is formalized as a particle process; we refer again to \cite{SW08}. 

Since we are considering only simple processes, it is convenient to identify such a process, which by definition is a counting measure, with its support, which is a locally finite set. In particular, a realization of $\widehat X$ is also considered as a set of hyperplanes, and a realization of $X$ is considered as a set of polytopes. The notations $\widehat X(\{H\})=1$ and $H\in \widehat X$ for a hyperplane $H$, for example, are therefore used synonymously.

The combinatorial type of a polytope $P$ in $\R^d$ is the set of all polytopes in $\R^d$ that are combinatorially isomorphic to $P$. Now we can formulate our result.

\begin{theorem}\label{T1}
Let $X$ be a tessellation of $\R^d$ that is induced by a stationary Poisson hyperplane process $\widehat X$ with spherical directional distribution $\varphi$. Suppose that the support of $\varphi$ is the whole unit sphere $\Sd$ and that $\varphi$ assigns measure zero to each great subsphere of $\Sd$. Then, with probability one, each combinatorial type of a simple $d$-polytope appears with positive density in $X$.
\end{theorem}

Theorem \ref{T1} implies, trivially, that under its assumptions almost surely each combinatorial type of a simple $d$-polytope appears infinitely often in $X$. When the latter fact was proved, among other results, in \cite{RS16}, a tool was a strengthened version of the Borel--Cantelli lemma, due to Erd\"os and R\'{e}nyi \cite{ER59} (see also \cite[p. 327]{Ren66}). When the note \cite{RS16} was submitted, an anonymous referee wrote ``that the use of ergodicity of the mosaic could lead to a possibly shorter alternative proof'', and he/she briefly indicated a possible approach. After thorough consideration, we preferred the more elementary Borel--Cantelli lemma. However, reconsideration revealed that ergodicity, applied in a different way, might lead to a stronger result, as far as the occurrence of combinatorial types is concerned. This is carried out in the following.

\section{Proof}\label{sec3}

Let $X$ satisfy the assumptions of Theorem \ref{T1}. Under the only assumption that the spherical directional distribution of the stationary Poisson hyperplane tessellation $X$ is zero on every great subsphere, it was shown in \cite[Thm. 10.5.3]{SW08} that $X$ is mixing and hence ergodic. This requires a few explanations. To model $X$ as a point process, we consider the space $\K$ of convex bodies (nonempty, compact, convex subsets) in $\R^d$ with the Hausdorff metric. By $\B(\K)$ we denote the $\sigma$-algebra of Borel sets in $\K$. Let ${\sf N}_s(\K)$ be the set of simple, locally finite counting measures on $\B(\K)$ and ${\mathcal N}_s(\K)$ its usual $\sigma$-algebra (for details see, e.g., \cite[Sect. 3.1]{SW08}). As underlying probability space $(\Omega,{\mathcal A},\bP)$, on which $X$ is defined, we can use $({\sf N}_s(\K),{\mathcal N}_s(\K),\bP_X)$, where $\bP_X$ is the distribution of $X$. For $t\in\R^d$, a bijective map ${\sf T}_t:\eta\mapsto {\sf T}_t\eta$ of ${\sf N}_s(\K)$ onto itself is defined by
$$ ({\sf T}_t\eta)(B):= \eta(B-t),\quad B\in \B(\K),\,\eta\in {\sf N}_s(\K).$$
Since $X$ is stationary, we have
$$ \bP_X({\sf T}_tA)=\bP_X(A)\quad\mbox{for }A\in{\mathcal N}_s(\K),$$
thus ${\sf T}_t$ induces a measure preserving map of ${\mathcal N}_s(\K)$ into itself. Let ${\mathcal T}:=\{ {\sf T}_t:t\in\R^d\}$. As shown in \cite[Thm. 10.5.3]{SW08}, the dynamical system $({\sf N}_s(\K),{\mathcal N}_s(\K),\bP_X,{\mathcal T})$ is mixing, that is,
$$ \lim_{\|t\|\to\infty} \bP_X(A\cap {\sf T}_tB)=\bP_X(A)\bP_X(B)$$
holds for all $A,B\in{\mathcal N}_s(\K)$. It follows that the system is ergodic, which means that $\bP_X(A)\in\{0,1\}$ for all $A\in {\bf T}:= \{A\in{\mathcal N}_s(\K):{\sf T}_tA=A\mbox{ for all }t\in\R^d\}$. Therefore, the `Individual Ergodic Theorem for $d$-dimensional Shifts' yields the following.

\begin{proposition}\label{P1}
Let $f$ be an integrable random variable on $({\sf N}_s(\K),{\mathcal N}_s(\K),\bP_X)$. Then
$$ \lim_{n\to\infty} \frac{1}{\lambda_d(B_n)} \int_{B_n} f({\sf T}_t\,\omega)\,\lambda(\D t) = \bE f$$
holds for $\bP_X$-almost all $\omega\in {\sf N}_s(\K)$.
\end{proposition}

We refer to Daley and Vere--Jones \cite[Proposition 12.2.II]{DV08} for a more general formulation (with hints to proofs of more general results in Tempel'man \cite{Tem72}). However, we have already incorpated into our Proposition \ref{P1} the information that in our case $({\sf N}_s(\K),{\mathcal N}_s(\K),\bP_X,{\mathcal T})$ is ergodic, which yields that the limit is equal to the expectation of $f$.

We apply this Proposition in the following way. First we choose a center function $c$ on $\K$; for example, let $c(K)$ denote the circumcenter of $K\in \K$, which is the center of the smallest ball containing $K$. Let $A\in \B(\K)$ be a translation invariant Borel set of convex bodies. Given any bounded Borel set $B\in\B(\R^d)$, we define
$$ f(B,\omega):= \sum_{K\in X(\omega),\,c(K)\in B} {\mathbbm 1}_A(K)$$
for $\omega\in\Omega$, where we use $(\Omega,{\mathcal A},\bP)=({\sf N}_s(\K),{\mathcal N}_s(\K),\bP_X)$ as the underlying probability space. Then $f(B,\cdot)$ is measurable, and $f(B+t,\omega)= f(B,{\sf T}_{-t}\,\omega)$ for $t\in\R^d$. The following generalizes an approach of Cowan \cite{Cow80} in the plane (``Tricks with small disks''). Assuming that $n>1$, we have

\begin{eqnarray*}
&& \int_{B_{n-1}} f(B_1+t,\omega)\,\lambda_d(\D t)\\
&&= \sum_{K\in X(\omega)} \int_{\R^d} {\mathbbm 1}\{t\in B_{n-1}\}{\mathbbm 1}\{K\in A\}{\mathbbm 1}\{c(K)\in B_1+t\}\,\lambda_d(\D t).
\end{eqnarray*}
Since
$$  {\mathbbm 1}\{t\in B_{n-1}\}{\mathbbm 1}\{c(K)\in B_1+t\}\le {\mathbbm 1}\{t\in -B_1+c(K)\}{\mathbbm 1}\{ c(K)\in B_n\},$$we get
\begin{eqnarray*}
&& \int_{B_{n-1}} f(B_1+t,\omega)\,\lambda_d(\D t)\\
&&\le \sum_{K\in X(\omega)} \int_{\R^d} {\mathbbm 1}\{t\in -B_1+c(K)\}{\mathbbm 1}\{K\in A\}{\mathbbm 1}\{c(K)\in B_n\}\,\lambda_d(\D t)\\
&& = \lambda_d(B_1)f(B_n,\omega).
\end{eqnarray*}
Similarly,
\begin{eqnarray*}
&& \int_{B_{n+1}} f(B_1+t,\omega)\,\lambda_d(\D t)\\
&& \ge \sum_{K\in X(\omega)} \int_{\R^d} {\mathbbm 1} \{t\in -B_1+c(K)\} {\mathbbm 1} \{K\in A\}{\mathbbm 1} \{c(K)\in B_n\}\,\lambda_d(\D t)\\
&& = \lambda_d(B_1)f(B_n,\omega).
\end{eqnarray*}
We conclude that
\begin{eqnarray*}
&& \frac{\lambda_d(B_{n-1})} {\lambda_d(B_n)} \frac{1}{\lambda_d(B_{n-1})} \int_{B_{n-1}} f(B_1,{\sf T}_{-t}\,\omega)\,\lambda_d(\D t)\\
&& \le \frac{\lambda_d(B_1)}{\lambda_d(B_n)}f(B_n,\omega)\\
&& \le \frac{\lambda_d(B_{n+1})}{\lambda_d(B_n)} \frac{1}{\lambda_d(B_{n+1})} \int_{B_{n+1}} f(B,{\sf T}_{-t}\,\omega)\,\lambda_d(\D t).
\end{eqnarray*}
By the Proposition, the lower and the upper bound converge, for $n\to\infty$, almost surely to $\bE f(B_1,\cdot)$, hence a.s.
\begin{equation}\label{3.1}
\lim_{n\to\infty} \frac{1}{\lambda_d(B_n)}f(B_n,\cdot) = \frac{\bE f(B_1,\cdot)}{\lambda_d(B_1)}.
\end{equation}

Now we assume in addition that there is a constant $D>0$ such that all convex bodies $K\in A$ satisfy ${\rm diam}\,K\le D$, where ${\rm diam}$ denotes the diameter. The center function $c$ satisfies $c(K)\in K$, hence if $c(K)\in B_{n-D}$ (with $n>D$) and ${\rm diam}\,K\le D$, then $K\subset B_n$. It follows that, for $n>D$,

\begin{eqnarray*}
&& \frac{\lambda_d(B_{n-D})} {\lambda_d(B_n)} \frac{1}{\lambda_d(B_{n-D})} \sum_{K\in X} {\mathbbm 1}_A(K) {\mathbbm 1}\{c(K)\in B_{n-D}) \\
&& \le\frac{1}{\lambda_d(B_n)} \sum_{K\in X,\,K\subset B_n} {\mathbbm 1}_A(K) \\
&&\le \frac{1}{\lambda_d(B_n)} \sum_{K\in X} {\mathbbm 1}_A(K\in A){\mathbbm 1}\{c(K)\in B_n\}.
\end{eqnarray*}

As $n\to\infty$, the lower and the upper bound converge a.s. to the right side of (\ref{3.1}), hence a.s. we have
\begin{equation}\label{3.1a} 
\delta(X,A):=\lim_{n\to\infty} \frac{1}{\lambda_d(B_n)} \sum_{K\in X,\,K\subset B_n} {\mathbbm 1}_A(K) = \frac{1}{\lambda_d(B^d)} \,\bE \sum_{K\in X,\,c(K)\in B^d} {\mathbbm 1}_A(K).
\end{equation}

Now we consider the special case where $A_D$ is the set of polytopes that are combinatorially isomorphic to a given simple $d$-polytope $P$ and have diameter at most $D$, for some fixed number $D>0$. We remark that (\ref{3.1a}) shows that
\begin{equation}\label{3.1b} 
\delta(X,A_D) = \frac{1}{\lambda_d(B^d)} \,\bE \sum_{K\in X,\,c(K)\in B^d} {\mathbbm 1}\{K\in A_D\},
\end{equation}
It remains to show that 
\begin{equation}\label{3.2}
\bE \sum_{K\in X,\,c(K)\in B^d} {\mathbbm 1}\{K\in A_D\}>0.
\end{equation}
For this, we use an argument from \cite{RS16}, which we recall for completeness. 

Without loss of generality, we can assume that $c(P) =o$. Let $F_1,\dots,F_m$ be the facets of $P$. We denote by $B(x,\varepsilon)$ the ball with center $x$ and radius $\varepsilon>0$, set $[B(x,\varepsilon)]_\Ha:= \{H\in\Ha:H\cap B(x,\varepsilon)\not=\emptyset\}$, and define
$$ A_j(P,\varepsilon) := \bigcap_{v\in{\rm vert}F_j} [B(v,\varepsilon)]_\Ha,\quad j=1,\dots,m,$$
where ${\rm vert}$ denotes the set of vertices. Each hyperplane from $A_j(P,\varepsilon)$ is said to be {\em $\varepsilon$-close} to $F_j$. A polytope $Q$ is said to be {\em $\varepsilon$-close} to $P$ if it has $m$ facets $G_1,\dots,G_m$ and, after suitable renumbering, the affine hull of $G_j$ is $\varepsilon$-close to $F_j$, for $j=1,\dots,m$. Since $P$ is simple and $c(P)=o$, we can choose numbers $D,\varepsilon_0>0$ such that for $0<\varepsilon\le \varepsilon_0$, the following is true:\\[1mm]
$\bullet$ the sets $A_1(P,\varepsilon),\dots,A_m(P,\varepsilon)$ are pairwise disjoint, and any hyperplanes $H_j\in A_j(P,\varepsilon)$, $j=1,\dots,m$, are the facet hyperplanes of a polytope $Q$ that is $\varepsilon$-close to $P$.\\[1mm]
$\bullet$ Any polytope $Q$ that is $\varepsilon$-close to $P$ satisfies the following:\\[1mm]
$\bullet$ $Q$ is combinatorially isomorphic to $P$,\\[1mm]
$\bullet$ $Q\subset P+B^d$,\\[1mm]
$\bullet$ ${\rm diam}\,Q\le D$,\\[1mm]
$\bullet$ $c(Q)\in B^d$.\\[1mm]
That this can be achieved by suitable choices of $D$ and $\varepsilon_0$, follows from easy continuity considerations and the fact that $P$ is simple. 

Now we define
$$ C(P,\varepsilon):= \{H \in \Ha: H\cap (P+B^d)\not=\emptyset,\; H\notin A_j(P,\varepsilon) \mbox{ for }j=1,\dots,m\}$$
and consider the event $E(P,\varepsilon)$ defined by
$$ \widehat X(A_j(P,\varepsilon))=1 \mbox{ for }j=1,\dots,m \enspace \mbox{and} \enspace\widehat X(C(P,\varepsilon))=0.$$
Let $0<\varepsilon\le\varepsilon_0$. The following was proved in \cite{RS16}:\\[1mm]
$\bullet$ If the event $E(P,\varepsilon)$ occurs, then some polytope $Q$ of the tessellation $X$ is $\varepsilon$-close to $P$ and hence satisfies $Q\in A_D$ and $c(Q)\in B^d$,\\[1mm]
$\bullet$  The event $\bP(E(P,\varepsilon))$ has positive probability.

Now it follows that 
$$ \bE \sum_{K\in X,\,c(K)\in B^d} {\mathbbm 1}\{K\in A_D\} \ge \bP(E(P,\varepsilon))>0,$$
which proves (\ref{3.2}).

The result is that $\delta(X,A_D)>0$ a.s. This implies, in particular, that with probability one the polytopes of the combinatorial type of $P$ appear in $X$ with positive density. Since there are only countably many combinatorial types, it also holds with probability one that each combinatorial type of a simple $d$-polytope appears in $X$ with positive density.

\noindent Author's address:\\[2mm]
Rolf Schneider\\
Mathematisches Institut, Albert-Ludwigs-Universit{\"a}t\\
D-79104 Freiburg i. Br., Germany\\
E-mail: rolf.schneider@math.uni-freiburg.de

\end{document}